\documentclass[a4paper,11pt]{article}
\usepackage{hyperref}
\hypersetup{
    colorlinks=true,
    linktoc=all,
    linkcolor=black,     
    citecolor=blue,     
}
\usepackage{graphicx}
\usepackage{amssymb}
\usepackage{latexsym, bm}
\usepackage{multicol}
\usepackage{indentfirst}
\usepackage{amsfonts}
\usepackage{amsmath}
\usepackage{setspace}
\newtheorem{theorem}{Theorem}[section]
\newtheorem{lemma}{Lemma}[section]
\newtheorem{claim}{Claim}[section]

\newtheorem{corollary}{Corollary}[section]

\newtheorem{problem}{Problem}[section]

\newcommand{\ignore}[1]{}

\begin{document}
\begin{spacing}{1}

\title{Sharp Ramsey thresholds for large books}
\date{}

\author{
Qizhong Lin,\footnote{Center for Discrete Mathematics, Fuzhou University,
Fuzhou, 350108, P.~R.~China. Email: {\tt linqizhong@fzu.edu.cn}. Supported in part  by National Key R\&D Program of China (Grant No. 2023YFA1010202), NSFC (No.\ 12171088, 12226401) and NSFFJ (No. 2022J02018).} \;\;\;\; Ye Wang\footnote{College of Mathematical Sciences, Harbin Engineering University, Harbin 150001, China. Email: {\tt ywang@hrbeu.edu.cn}. Supported in part by NSFC (No. 12101156).}
}

\maketitle

\begin{abstract}

For graphs $G$ and $H$, let $G\to H$ signify that any red/blue edge coloring of $G$ contains a monochromatic $H$. Let $G(N,p)$ be the random graph of order $N$ and edge probability $p$. The Ramsey thresholds for fixed graphs have received most attention.
In this paper, we consider the Ramsey thresholds in another angle. In particular, we will consider the sharp Ramsey threshold for the large book graph $B_n^{(k)}$, which consists of $n$ copies of $K_{k+1}$ all sharing a common $K_k$. In particular, for every fixed integer $k\ge 2$ and for any real $c>1$, let $N=c2^k n$. Then for any real $\gamma>0$,
\[
\lim_{n\to \infty} \Pr(G(N,p)\to B_n^{(k)})= \left\{  \begin{array}{cl}
0  & \mbox{if $p\le\frac{1}{c^{1/k}}(1-\gamma)$,} \\  1 & \mbox{if $p\ge\frac{1}{c^{1/k}}(1+\gamma)$}. \end{array}      \right.
\]
This implies that $r(B_n^{(k)},B_n^{(k)})=2^kn+o(n)$, and hence especially extends the work of Conlon (2019) and the follow-up work of Conlon, Fox and Wigderson (2022) on book Ramsey numbers.

\medskip

\end{abstract}

\section{Introduction}

For graphs $G$ and $H$, let $G\to H$ signify that any red/blue edge coloring of $G$ contains a monochromatic copy of $H$. The Ramsey number $r(H)$ is defined as the minimum $N$ such that $K_N\to H$.
Ramsey's theory \cite{ram} guarantees that the Ramsey number $r(H)$ is finite for all $H$.
The question of whether or not $G$ has the Ramsey property $G\to H$ is of particular interest when $G$ is a typical random graph from the probability space ${\cal G}(n,p)$, defined by Erd\H{o}s-R\'enyi \cite{erd-ren}, where $n$ is the number of ordered vertices and $p$ is the probability of edge appearance. A random graph in ${\cal G}(n,p)$ is always denoted by $G(n,p)$.

The Ramsey threshold $p(n)$ of the event $G(n,p)\to H$ is defined by
\[
\lim_{n\to \infty} \Pr(G(n,p)\to H)= \left\{  \begin{array}{cl}
0  & \mbox{if $p\ll p(n)$,} \\  1 & \mbox{if $p\gg p(n)$}. \end{array}      \right.
\]
We also call $p_\ell=o(p(n))$ and $p_u=\Omega(p(n))$ a lower Ramsey threshold and an upper Ramsey threshold, respectively.
It is often to signify  $\Pr(G(n,p)\to H)\to 1$ and   $\Pr(G(n,p)\to H)\to 0$ as $n\to\infty$
by saying that {\em asymptotically almost surely} (a.a.s.) $G(n,p)\to H$ and a.a.s. $G(n,p)\not\to H$, respectively.
If we can replace $p\ll p(n)$ and $p\gg p(n)$ in the above with $p\le (1-\gamma)p(n)$ and $p\ge (1+\gamma)p(n)$ for every $\gamma>0$, respectively, then the Ramsey threshold is said to be {\bf sharp}. For convenience, we always say such $p(n)$ is a sharp Ramsey threshold for $H$.

The study of Ramsey thresholds was initiated by Frankl and R\"{o}dl \cite{f-r} and independently by {\L}uczak, Ruci\'nski and Voigt \cite{l-r-v}, who proved that $p =1/\sqrt{n}$ is a Ramsey threshold for triangle.
In a series of papers \cite{f-r,l-r-v,r-r93,r-r-1,r-r}, the Ramsey thresholds are determined for any fixed graph $H$.
For a graph $H$, let $v(H)$ and $e(H)$ be the numbers of vertices and edges of $H$, respectively.
The Ramsey threshold for a fixed graph was determined by R\"{o}dl and Ruci\'nski \cite{r-r}, who proved that (except $H$ is a path of length 3 as was pointed out in \cite{f-k} or a disjoint union of stars)
\begin{align}\label{r-r}
\lim_{n\to \infty} \Pr(G(n,p)\to H)= \left\{  \begin{array}{cl}
0  & \mbox{if $p\ll n^{-1/m_2(H)}$,} \\  1 & \mbox{if $p\gg n^{-1/m_2(H)}$}, \end{array}      \right.
\end{align}
where $m_2(H)=\max\{\frac{e(F)-1}{v(F)-2}:\; F\subseteq H,\; v(F)\ge 3\}$. This result has a short proof from Nenadov and Steger \cite{n-s}.

The sharp thresholds for Ramsey properties seemed out of hand until a general technique for settling these questions was introduced by Friedgut \cite{f99}. In particular, Friedgut and Krivelevich \cite{f-k} obtained all sharp thresholds for fixed trees except the star and a path of length 3. When $H$ is a triangle, it was established by Friedgut, R\"{o}dl, Ruci\'{n}ski and Tetali \cite{frrt}.

In the following, we mainly focus on the situations when the graphs are large.
A closely related problem is the size Ramsey number.
For a graph $H$, Erd\H{o}s, Faudree, Rousseau and Schelp \cite{efrs-size} defined the size Ramsey number as $\hat{r}(H)=\min\{e(G): G\to H\}$.

Beck \cite{beck-83} proved $\hat{r}(P_n)=O(n)$ for path $P_n$ of length $n$, who in fact showed a.a.s. $$G(c_1n,c_2/n)\to P_n,$$ where $c_1$ and $c_2$ are positive constants. This has been improved by Dudek and Pra{\l}at \cite{d-p}.
For a path $P_n$ of length $n$,  Gerencs\'er and Gy\'arf\'as \cite{g-g} proved
$
r(P_n)=n+\lceil\frac{n}{2}\rceil.
$
Thus  if $N< 3n/2$, $G(N,1)\to P_n$ is an impossible event.
Letzter \cite{let} proved that if $c>1$, then a.a.s. $$G(3cn/2,p)\to P_n$$ provided $pn\to \infty$, hence $p=\frac{1}{n}$ is a Ramsey threshold of $G(3cn/2,p)\to P_n$, which improves Beck's result \cite{beck-83} further.

Let ${\cal F}_{\Delta,n}$ be the family of graphs $H$ with order $n$ and maximum degrees at most $\Delta$. 
Beck \cite{beck-2} conjectured that the size Ramsey number $\hat{r}(H)=O(n)$ for any $H\in {\cal F}_{\Delta,n}$.
However, R\"{o}dl and Szemer\'edi \cite{r-s} showed that it does not hold even for $\Delta=3$.
In 2011, Kohayakawa,  R\"{o}dl, Schacht and Szemer\'edi \cite{k-r-s-s} proved that for every fixed $\Delta\ge2$, there exist constants $B=B(\Delta)$ and $C=C(\Delta)$ such that if $N= \lceil Bn\rceil$ and $p=C(\log N/N)^{1/\Delta}$, then for any $H\in{\cal F}_{\Delta,n}$,
\[
\lim_{n\to\infty}\Pr\left(G(N,p)\to H\right)=1.
\]
This implies that $\hat{r}(H)=O(n^{2-1/\Delta}\log^{1/\Delta}n)$ for any $H\in{\cal F}_{\Delta,n}$.

Let $B_{n}^{(k)}$ be the book graph consisting of $n$ copies of $K_{k+1}$, all sharing a common $K_k$.
Let $K_{k,n}$ be the complete bipartite graph with two parts of sizes $k$ and $n$. Clearly, both of these two special families of graphs do not belong to ${\cal F}_{\Delta,n}$. The corresponding Ramsey-type problems of these two families of graphs have attracted a great deal of attention.
Li, Tang and Zang \cite{ltz} proved that for every fixed $k\ge2$, $r(K_{k,n})=(2^k+o(1)) n$, but $r(B_{n}^{(k)})$ is much harder to handle.
Erd\H{o}s, Faudree, Rousseau and Schelp \cite{efrs-size} and independently Thomason \cite{tho} proved that $(2^k+o(1))n\le r(B_n^{(k)},B_n^{(k)})\le4^kn.$
After many years, Conlon \cite{Conlon} proved that for every $k\ge2$,
\begin{align}\label{con-bk}
r(B_{n}^{(k)})=(2^k+o(1)) n,
\end{align}
which confirms a conjecture of Thomason \cite{tho} asymptotically and also gives an answer to a problem proposed by Erd\H{o}s \cite{efrs-size}.
The small term $o(1)$ in (\ref{con-bk}) is improved further by Conlon, Fox and Wigderson \cite{cfw} using a different method.
Books have attracted a great
deal of attention in graph Ramsey theory, see e.g. the recent breakthrough of Campos,
Griffiths, Morris and Sahasrabudhe \cite{cgms}.
For more Ramsey numbers of books, the reader is referred to \cite{cl,cly,cfw-2,fly,n1,n2,n3,rou-sh} etc.

A classical result of Erd\H{o}s, Faudree, Rousseau and Schelp \cite{efrs-size} also tells that
$\hat{r}(K_{k,n})=\Theta(n)$, and $\hat{r}(B_{n}^{(k)})=\Theta(n^2).$
 Recently, Conlon, Fox and Wigderson \cite{cfw-s} establish that  for every fixed $k\ge2$ and all large $n$,
 $$\hat{r}(B_{n}^{(k)})=\Theta(k2^kn^2).$$
 Moreover, they obtain that if $n=\Omega(k\log k)$, then $\hat{r}(K_{k,n})=\Theta(k^22^kn).$

In this paper, we mainly obtain sharp Ramsey thresholds for $B_{n}^{(k)}$ and $K_{k,n}$.

\begin{theorem}\label{main}
Let $N=c2^k n$, where  $k\ge 2$ is an integer and $c>1$ is a real number. Then for any $\gamma>0$,
\[
\lim_{n\to \infty} \Pr(G(N,p)\to B_n^{(k)})= \left\{  \begin{array}{cl}
0  & \mbox{if $p\le\frac{1}{c^{1/k}}(1-\gamma)$,} \\  1 & \mbox{if $p\ge\frac{1}{c^{1/k}}(1+\gamma)$}. \end{array}      \right.
\]
\end{theorem}

{\em Remark.}
Theorem \ref{main} extends (\ref{con-bk}). Indeed, if we take $c=1+\epsilon$ for sufficiently small $\epsilon>0$ and $p\to1$, then a.a.s. $G((1+\epsilon)2^kn,p)\to B_n^{(k)}$. Especially, (\ref{con-bk}) holds.
Moreover, Theorem \ref{main} implies that the sharp Ramsey threshold for $B_n^{(k)}$ is a positive constant $1/c^{1/k}$, although the edge density of the book graph $B_n^{(k)}$ tends to zero as $n\to \infty$.

\ignore{
For $c>1$, since the Ramsey threshold of $G(cr(P_n),p)\to P_n$ is $\frac{1}{n}$ in Letzter \cite{let}
and the sharp Ramsey threshold  of $G(cr(K_{1,n}),p)\to K_{1,n}$ is $\frac{1}{c}$ from Theorem \ref{main},
we see that the Ramsey thresholds of $G(cr(T_n),p)\to T_n$ are so different for different types of trees $T_n$.
Note that most Ramsey numbers $r(T_n)$ of trees $T_n$ with $n$ edges are unknown, ranging from $\lfloor\frac{4n+1}{3}\rfloor$ to $2n$, see Erd\H{o}s, Faudree, Rousseau and Schelp \cite{efrs-b}, and Yu and Li \cite{y-l}.
It would be interesting to find more Ramsey thresholds for different types of trees. Moreover, it would be very interesting to determine the sharp Ramsey threshold for $P_n$, which seems not easy.
}

\medskip
Combining Lemma \ref{main-le}, the following is immediate.
\begin{corollary}
Let $N=c2^k n$, where  $k\ge 2$ is an integer and $c>1$ is a real number. Then for any $\gamma>0$,
\[
\lim_{n\to \infty} \Pr(G(N,p)\to K_{k,n})= \left\{  \begin{array}{cl}
0  & \mbox{if $p\le\frac{1}{c^{1/k}}(1-\gamma)$,} \\  1 & \mbox{if $p\ge\frac{1}{c^{1/k}}(1+\gamma)$}. \end{array}      \right.
\]
\end{corollary}

\medskip\noindent
{\bf  Notation:} For a graph $G=(V,E)$ with vertex set $V$ and edge set $E$, let $uv$ denote an edge of $G$. For $X \subseteq V$, $e(X)$ is the number of edges in $X$, and $G[X]$ denotes the subgraph of $G$ induced by $X$. For two disjoint subsets $X,Y\subseteq V$, $e_G(X,Y)$ denotes the number of edges between $X$ and $Y$. In particular, the neighborhood of a vertex $v$ in $U\subseteq V$ is denoted by $N_G(v,U)$, and $\deg_G(v,U)=|N_G(v,U)|$ and the degree of $v$ in $G$ is $\deg_G(v)=|N_G(v,V)|$.
Let $X \sqcup Y$ denote the disjoint union of $X$ and $Y$. We always delete the subscriptions if there is no confusion from the context. Note that we have not distinguished large $x$ from $\lceil x\rceil$ or $\lfloor x\rfloor$ when $x$ is supposed to be an integer since these rounding errors are negligible to the asymptotic calculations.

\medskip
The rest of the paper is organized as follows. In Section \ref{low-b}, we will give the proof of the lower Ramsey threshold  of Theorem \ref{main}. In Section \ref{up-b}, we shall present the proof of the upper Ramsey threshold of Theorem \ref{main}.  Finally, we will have some discussions in Section \ref{clu}.

\section{The lower Ramsey threshold}\label{low-b}

We need the Chernoff's bound in the following form, see e.g. \cite{alo-spe,Bollobas,cher,jlr,li-l}.

\begin{lemma}\label{cher}
Let $X_1,X_2,\dots,X_n$ be mutually independent variables such that $\Pr(X_i=1)=p$ and $\Pr(X_i=0)=1-p$ for $1\le i\le n$
and $S_n=\sum_{i=1}^n X_i$. Then, there exists $\delta_0 >0$ such that
\[
\Pr\big[S_n \ge n(p+\delta)\big] < \exp\Big( -\frac{n\delta^2}{3pq}\Big)
\]
for any $\delta \in (0,\delta_0)$.
\end{lemma}

The following slightly stronger lemma implies the lower Ramsey threshold of Theorem \ref{main}.

\begin{lemma}\label{main-le}
Let $k\ge 2$ be an integer and $\gamma\in (0,1)$. Let $c=c(n)$  be a function
such that $1< c(n)\le e^{o(n)}$ and $N=c 2^k n$. If $p\le\frac{1}{c^{1/k}}(1-\gamma)$, then
\[
\lim_{n\to \infty} \Pr\big(G(N,p)\to K_{k,n}\big) = 0.
\]
\end{lemma}
{\bf Proof.} It suffices to show that
\[
p_{_\ell}=\frac{1}{c^{1/k}}(1-\gamma)^{1/k}
\]
is a function such that a.a.s. $G(N,p_{_\ell)}\not\to K_{k,n}$.
\begin{claim}\label{half}
Let $p_0=\frac{p_{_\ell}}{2}$. Then a.a.s. $G(N,p_0)$ contains no $K_{k,n}$.
\end{claim}
{\bf Proof of Claim \ref{half}.}  Let $V$ be  the vertex set with $|V|=N$. Consider the random graph $G(N,p_0)$ on $V$.
Let $U\subseteq V$ be a subset with $|U|=k$ and $V\setminus U=\{v_1,v_2,\dots,v_{N-k}\}$.
For $1\le i\le N-k$, define a random variable $X_i$ such that $X_i=1$
if $v_i$ is a common neighbor of $U$ and $0$ otherwise. Then
\[
\Pr(X_i=1)=p_0^k = \frac{1}{c 2^k}(1-\gamma)
\]
and $\Pr(X_i=0)=1-p_0^k$.

Set a random variable $S_{N-k}=\sum _{i=1}^{N-k}X_i$ that has the binomial distribution $B(N-k,p_0^k)$.
Note that the event $S_{N-k}\ge n$ means that $G(N,p_0)$ contains $K_{k,\,n}$ with $U$ as the part of $k$ vertices.
Hence
\[
\Pr\big(K_{k,n}\subseteq G(N,p_0)\big) \le \binom{N}{k}\Pr(S_{N-k}\ge n).
\]
To evaluate $\Pr(S_{N-k}\ge n)$, let us write
the event $S_{N-k}\ge n$ as $S_{N-k}\ge (N-k)(p_0^k +\delta)$, where $n=\frac{N}{c 2^k} = (p_0^k +\delta)(N-k)$ and
\[
\delta =\frac{N}{c 2^k(N-k)}-p_0^k =\frac{N}{c 2^k(N-k)}-\frac{1}{c 2^k}(1-\gamma)= \frac{1}{c 2^k}\Big(\gamma+\frac{k}{N-k}\Big).
\]
Therefore, Lemma \ref{cher} implies
\[
\Pr(S_{N-k}\ge n) =\Pr\big(S_{N-k}\ge (p_0^k +\delta)(N-k)\big)\le \exp\left\{-\frac{(N-k)\delta^2}{3p_0^k(1-p_0^k)}\right\}.
\]
Note that $\delta\sim\frac{\gamma}{c 2^k}$ and
\[
(N-k)\delta^2 \sim  N \delta^2 = (c 2^k n) \delta^2 \sim \frac{\gamma^2 n}{c 2^k}.
\]
Hence we have
\[
\binom{N}{k}\Pr(S_{N-k}\ge n) \le N^k\exp\left\{-\frac{\gamma^2 n}{4 c 2^k p_0^k(1-p_0^k)}\right\}.
\]
Note that $c2^k p_0^k =1-\gamma$ and thus
\[
\frac{\gamma^2 n}{4 c 2^k p_0^k(1-p_0^k)}  =  \frac{\gamma^2 n}{4(1-\gamma)(1-p_0^k)}\ge \frac{\gamma^2 n}{4}.
\]
Since $c\le e^{o(n)}$, we have
\[
N^k =(c2^kn)^k = \exp\big[k(\log c +k\log 2 +\log n)\big] = e^{o(n)},
\]
hence
\[
\binom{N}{k}\Pr(S_{N-k}\ge n) \le \exp\left(-\frac{\gamma^2 n}{4} +o(n)\right) \to 0,
\]
and the claim follows.          \hfill$\Box$

\medskip

To finish the proof, we shall show that a.a.s. $G(N,p)\not\to K_{k,n}$.

Let us write the defined random variable $S_{N-k}$ as $S_{N-k}^{p_{_\ell}/2}(U)$ for fixed $U$ with $|U|=k$,
where the superscript $p_{_\ell}/2$ corresponds to random graph $G(N,p_{_\ell}/2)$.
Then we have shown
\begin{equation}\label{half-p}
\binom{N}{k}\Pr\left(S_{N-k}^{p_{_\ell}/2}(U)\ge n\right)  \to  0,
\end{equation}
as $n\to\infty$.
Consider an edge coloring of $G(N,p_\ell)$ with red and blue at random with probability $1/2$, independently.
It is easy to see that both red graphs and blue graphs form ${\cal G}(N,p_{_\ell}/2)$.

For a vertex set $U$ of size $k$, let $S_{N-k}^{p_{_\ell},R}(U)$ and $S_{N-k}^{p_{_\ell},B}(U)$ be the numbers of common red and blue neighbors of $U$, respectively.
Then
\[
\Pr\left(S_{N-k}^{p_{_\ell},R}(U)\ge n\right) = \Pr\left(S_{N-k}^{p_{_\ell},B}(U)\ge n\right)= \Pr\left(S_{N-k}^{p_{_\ell}/2}\ge n\right),
\]
and thus $\Pr\big[S_{N-k}^{p_{_\ell},R}(U)\ge n\;\;\mbox{or}\;\;S_{N-k}^{p_{_\ell},B}(U)\ge n\big] \le 2 \Pr\big[S_{N-k}^{p_{_\ell}/2}\ge n\big]$.
Therefore, from (\ref{half-p}), we have
\[
\binom{N}{k} \Pr\left(S_{N-k}^{p_{_\ell},R}(U)\ge n\;\;\mbox{or}\;\;S_{N-k}^{p_{_\ell},B}(U)\ge n \right) \to 0
\]
as $n\to\infty$, which implies that a.a.s. $G(N,p_\ell)\not\to K_{k,n}$.  \hfill $\Box$

\medskip

\section{The upper Ramsey threshold}\label{up-b}

The following result follows from Chernoff bound directly.
\begin{lemma}\label{quasi}
Let $p\in (0,1]$ be a fixed probability.
If $N\to\infty$, then a.a.s. $G\in {\cal G}(N,p)$ with vertex set $V$ satisfies the following properties:

\medskip
(i)  For any vertex $v\in V$ and subset $U\subseteq V$, $\deg(v,U)= p|U|+o(N)$;

(ii) For any pair of distinct vertices $u$ and $v$, $|N(u)\cap N(v)|= p^2 N+o(N)$;

(iii) For any subsets $U\subseteq V$, $e(U)=p{\binom{|U|}2}+o(N^2)$;

(vi) For any disjoint vertex sets $U$ and $W$, $e(U,W)=p|U||W|+o(N^2)$.
\end{lemma}


\subsection{The first case for $k=2$}

In this subsection, we include a short proof for the case when $k=2$ of Theorem \ref{main}.
Denote $B_n$ instead of $B_n^{(2)}$. The upper Ramsey threshold for $k=2$  follows from the following lemma.

\begin{lemma}\label{le-2b}
Let $c_0>1$ be a constant. Let $c=c(n)\ge c_0$ and $p=\frac{1+\gamma}{\sqrt{c}}$, where $\gamma\in(0,\sqrt{c_0}-1]$.
If $G$ is a graph of order $N= 4cn$ that satisfies properties in Lemma \ref{quasi}, then  $G\to B_n$ for all large $n$.
\end{lemma}
{\bf Proof.} Suppose that there is an edge-coloring of $G$ by red and blue that contains no monochromatic $B_n$.
We shall show this assumption would lead to a contradiction.

Let $V$ be the vertex set of $G$. Let $R$ and $B$ denote the red and blue subgraphs, respectively.
Let $M_r$ and $M_b$ be the number of monochromatic triangles in red and blue, respectively.
Let $M_{rb}$  be the numbers of non-monochromatic triangles. Denote by $M=M_r+M_b$ the number of  monochromatic triangles,
and $T=M+M_{rb}=M_r+M_b+M_{rb}$ the number of triangles in $G$.

 Note from Lemma \ref{quasi} that $e(G)\sim\frac{1}{2}p N^2$,
and $|N(u)\cap N(v)|\sim p^2 N$, we have
\begin{equation}\label{tri}
T=\frac{1}{3}\sum_{uv\in E(G)} |N(u)\cap N(v)| \sim \frac{1}{6} p^3 N^3,
\end{equation}
where coefficient $\frac{1}{3}$ of the sum follows from that each triangle is counted triply in the sum.

Since a red edge $uv$ and $n$ red common neighbors of $u$ and $v$  yield a red $B_n$,
we have $|N_R(u)\cap N_R(v)| \le n-1$. Hence
\[
M_r = \frac{1}{3} \sum_{uv\in E(R)} |N_R(u)\cap N_R(v)| \le \frac{1}{3}(n-1)e(R).
\]
Similarly, $M_b \le \frac{1}{3}(n-1)e(B)$, and thus
\begin{equation}\label{mono}
M \le \frac{1}{3}(n-1)e(G) \sim \frac{1}{6} p n N^2.
\end{equation}
For any $v\in V$, each edge between $N_R(v)$ and $N_B(v)$ is contained in a non-monochromatic triangle, and thus
\[
M_{rb} = \frac{1}{2} \sum_{v\in V} e(N_R(v),N_B(v)) =  \frac{1}{2} \sum_{v\in V} p \deg_R(v)\deg_B(v) +o(N^3),
\]
where $\frac{1}{2}$ comes from that each such triangle is counted by its two vertices and the term $o(N^3)$ comes from the third property in Lemma \ref{quasi}.
Since $\deg_R(v)+\deg_B(v)=\deg(v)$, we have $\deg_R(v)\deg_B(v)\le \frac{1}{4} [\deg(v)]^2$. Therefore,
\begin{equation}\label{non-m}
M_{rb} \le \frac{1}{8}p \sum_{v\in V} [\deg(v)]^2 +o(N^3)\sim \frac{1}{8} p^3 N^3.
\end{equation}
Recall $M=T-M_{rb}$, which and (\ref{tri}), (\ref{mono}) and (\ref{non-m}) yield
\[
\frac{1}{6} p n N^2\ge (1-o(1)) \left(\frac{1}{6} p^3 N^3 - \frac{1}{8} p^3 N^3 \right)=\left(\frac{1}{24}-o(1)\right) p^3 N^3,
\]
which implies that $p^2\le (1+o(1))\frac{4n}{N}=(1+o(1))\frac{1}{c}$,
contradicting to the assumption $p=\frac{1+\gamma}{\sqrt{c}}$ with $\gamma>0$ fixed,
and the proof is completed.  \hfill $\Box$

\medskip

\begin{theorem}\label{th-2book}
Let $c_0>1$ be a constant and $c=c(n)$ be a function such that $c_0\le c \le e^{o(n)}$. If $N=4cn$ and $\gamma\in (0,\sqrt{c_0}-1]$, then
\[
\lim_{n\to \infty} \Pr(G(N,p)\to B_n)= \left\{  \begin{array}{cl}
0  & \mbox{if $p\le\frac{1}{\sqrt{c}}(1-\gamma)$,} \\  1 & \mbox{if $p\ge\frac{1}{\sqrt{c}}(1+\gamma)$}. \end{array}      \right.
\]
\end{theorem}

The proof of Theorem \ref{th-2book} comes from Lemma \ref{main-le} and Lemma \ref{le-2b} immediately,
which yields a corollary as follows.

\begin{corollary}\label{cor-b}
If $n\ll N \le n e^{o(n)}$,  then the sharp Ramsey threshold of $B_n$ in $G(N,p)$ is $2\sqrt{\frac{n}{N}}$ as $n\to\infty$.
\end{corollary}

Let us mention that  from (\ref{r-r}), the threshold of $B_m$ for fixed $m$ in $G(N,p)$ is $\Theta{(1/\sqrt{N})}$ as $m_2(B_m)=2$.
It is natural to ask the following problem.

\begin{problem}
Prove or disprove that $2\sqrt\frac{m}{N}$ is the sharp threshold of $B_m$ in $G(N,p)$  for fixed $m$.
If affirmative, what about probability $p=2\big(\frac{m}{N}\big)^{1/k}$ and book $B_m^{(k)}$ for fixed $k\ge2$ and $m\ge 1$?
\end{problem}

In the following, we shall focus on the sharp threshold of $B_n^{(k)}$ in $G(N,p)$  for fixed $k\ge3$ and sufficiently large $n\ge 1$.

\subsection{The regularity method and useful lemmas}\label{sec:regularitytools}


Szemer\'{e}di regularity lemma \cite{szemeredi,szemeredi1} is a powerful tool in extremal graph theory. There are many important applications of the regularity lemma. We refer the reader to nice surveys \cite{kss,ko-sim,rs} and other related references. The proof for the upper Ramsey thresholds of Theorem \ref{main} for general $k\ge3$ mainly relies on the regularity method.

Given  $p \in (0,1]$ and $\varepsilon > 0$, the $p$-density of a pair $(U,W)$ of sets of vertices in a graph $G$ is defined as
$d_{G,p}(U,W) = \dfrac{e_G(U,W)}{p|U||W|}.$
We say that the pair $(U,W)$ is \emph{$(\varepsilon,p)$-regular} in $G$ if $|d_{G,p}(U,W) - d_{G,p}(U',W')| \le \varepsilon$ for all $U' \subset U$ and $W' \subset W$ with $|U'| \ge \varepsilon |U|$ and $|W'| \ge \varepsilon |W|$.
 When $p=1$, it is the usual edge density, denoted by $d_{G}(U,W)$, between $U$ and $W$. Also, the set $U$ is said to be $(\varepsilon,p)$-regular if the pair $(U,U)$ is $(\varepsilon,p)$-regular. We write $d(U)$ for $d(U,U)$.

Given $0<\eta,p\le1$, $D\ge1$, a graph $G$ is called $(\eta, p, D)$-upper-uniform if, for all disjoint
sets of vertices $U, W$ of size at least $\eta|V(G)|$, the density $d_{G,p}(U,W)$ is at most $D$.
Given a red-blue coloring of the edges of $G$, we write $R$ and $B$ for the graphs on $V(G)$ induced by the red and blue edges, respectively.
We say that $V(G) = \sqcup_{i=1}^m V_i$ is an \emph{equitable partition} for the coloring $(R,B)$ of $G$ if $\big| |V_i|-|V_j|\big|\le1$ for all $1\le i<j\le m$.

We will use the following regularity lemma for random graphs.

\begin{lemma}\label{reglem}
For any $\varepsilon > 0$ and integer $M_0\ge1$, there exists $M=M(\varepsilon,M_0)>M_0$ such that the following holds. If $p\in (0,1]$ is fixed, then {\bf a.a.s.} every $2$-coloring of the edges of $G\in\mathcal{G}(N,p)$ has an $(\varepsilon, p)$-regular equitable partition $V(G)=\sqcup_{i=1}^m V_i$ where $M_0 \leq m \leq M$ such that

\medskip
(i)	 each part $V_i$ is $(\varepsilon, p)$-regular;

(ii) for each $V_i$, all but at most $\varepsilon m$ parts $V_j$ such that $(V_i,V_j)$ are $(\varepsilon, p)$-regular;

(iii) for any vertex $v\in V$ and for $1\le i\le m$, $\deg(v,V_i)= p|V_i|+o(N)$;

(iv)  for $1\le i\le m$, $e(V_i)=p{\binom{|V_i|}2}+o(N^2)$;

(v)  for $1\le i<j\le m$, $e(V_i,V_j)=p|V_i||V_j|+o(N^2)$.

\end{lemma}
{\bf Proof.} We only sketch the proof of Lemma \ref{reglem} as follows. From Lemma \ref{quasi}, a.a.s. $G\in {\cal G}(N,p)$ satisfies that (1) for any vertex $v\in V$ and subset $U\subseteq V$, $\deg(v,U)= p|U|+o(N)$;
(2) for any subsets $U\subseteq V$, $e(U)=p{\binom{|U|}2}+o(N^2)$;
(3) for any disjoint vertex sets $U$ and $W$, $e(U,W)=p|U||W|+o(N^2)$. Therefore, the random graph $G$ and hence the red subgraph $R$ and the blue subgraph $B$ are a.a.s. upper uniform (with suitable parameters). Let $\varepsilon_1=\varepsilon/2$, $\varepsilon_2=\varepsilon^2/128$, $K_1=K(\varepsilon_1)\le2^{(1/\varepsilon_1)^{(10/\varepsilon_1)^{15}}}$, and let $\eta=\min\{\varepsilon_1/K_1,\varepsilon^3/256\}$ as in Conlon, Fox and Wigderson \cite[Lemma 2.1]{cfw}.
We can first apply the colored version of Letzter \cite[Theorem 5.2]{let} (from an original version by Kohayakawa and R\"{o}dl \cite{k-r1,k-r-2}) to obtain that there exists $L=L(\eta,M_0)>M_0$ such that the following holds. If $p\in (0,1]$ is fixed, then we have that a.a.s. every $2$-coloring of the edges of $G\in\mathcal{G}(N,p)$ has an equitable partition $V(G)=\sqcup_{i=1}^\ell W_i$ with $\max\{M_0,1/\eta\}\le \ell\le L$ such that all but at most $\varepsilon {\binom m2}$ pairs $(W_i,W_j)$ are $(\eta, p)$-regular. Then we apply \cite[Lemma 2.4]{cfw} to each $W_i$ to get an equitable partition
$W_i = U_{i1}\sqcup\cdots\sqcup U_{iK_1}$ such that each $U_{ij}$ for $1\le j\le K_1$ is $\varepsilon_1$-regular.
Subsequently, by a similar argument as that in \cite[Lemma 2.1]{cfw}, we can obtain an $(\varepsilon, p)$-regular equitable partition $V(G)=\sqcup_{i=1}^m V_i$ satisfying the conditions from the above equitable partition as desired.\hfill$\Box$

\medskip

The following is a standard counting lemma, see e.g. \cite[Theorem 18]{rs}.

\begin{lemma}[R\"{o}dl and Schacht \cite{rs}]\label{cl-RS}
For any $\eta>0$, there exists $\varepsilon>0$ such that if $V_1,\ldots,V_k$ are distinct subsets of a graph $G$ such that all pairs $(V_i,V_j)$ are $\varepsilon$-regular.
Then the number of labeled copies of $K_k$ whose $i$th vertex is in $V_i$ for all $i$ is at least
\[
\left( \prod_{1 \leq i<j\leq k}d(V_i,V_j) - \eta \right) \prod_{i=1}^k |V_i|.
\]
\end{lemma}

As a simple corollary, we have the following result.
\begin{lemma}\label{countinglemma}
For any $\eta>0$, there exists $\varepsilon>0$ such that if $V_1,\ldots,V_k$ are (not necessarily distinct) subsets of a graph $G$ such that all pairs $(V_i,V_j)$ are $\varepsilon$-regular.
Then the number of labeled copies of $K_k$ whose $i$th vertex is in $V_i$ for all $i$ is at least
\[
\left( \prod_{1 \leq i<j\leq k}d(V_i,V_j) - \eta \right) \prod_{i=1}^k |V_i|.
\]
\end{lemma}
{\bf Proof.} Suppose first that $V_1,\ldots,V_k$ are all the same, i.e., $V_1=V_2=\cdots=V_k$. Let $V_1=\{v_1,\dots,v_t\}$. Then we would set new distinct subsets $U_i:=\{v_1^i,\dots,v_t^i\}$ for $i\in[k]$ such that for $i\neq i'$ and $j\neq j'$, $v_j^i$ is adjacent with $v_{j'}^{i'}$ if and only if $v_jv_{j'}$ is an edge in $V_1$. From the assumption, we know that for all pairs $(U_j,U_{j'})$ are $\varepsilon$-regular for  $j\neq j'$ since $V_1$ is $\varepsilon$-regular. Therefore, from Lemma \ref{cl-RS}, the number of labeled copies of $K_k$ whose $i$th vertex is in $U_i$ for all $i$ is at least
$( \prod_{1 \leq i<j\leq k}d(U_i,U_j) - \eta) \prod_{i=1}^k |U_i|.$ From the definition of $U_i$, we know that a labeled copy of $K_k$ whose $i$th vertex is in $U_i$ for all $i$ also forms a labeled copy of $K_k$ in $V_1$.
Therefore, the number of labeled copies of $K_k$ in $V_1$ is at least
$( \prod_{1 \leq i<j\leq k}d(V_1) - \eta) \prod_{i=1}^k |V_1|.$ For all other cases, the arguments are similar.\hfill$\Box$

\medskip
The following is a counting lemma by Conlon \cite[Lemma 5]{Conlon}, which will be used to find a large monochromatic book.
\begin{lemma}[Conlon \cite{Conlon}]\label{con-ct}
For any $\delta>0$ and any integer $k\ge1$, there is $\varepsilon > 0$ such that if $V_1,\ldots,V_k$, $V_{k+1},\ldots,V_{k+\ell}$, are (not necessarily distinct) vertex with $(V_i,V_{i'})$ $\varepsilon$-regular  of density $d_{i,i'}$ for all $1\le i<i'\le k$ and $1\le i\le k<i'\le k+\ell$ and $d_{i,i'}\ge\delta$ for all $1\le i<i'\le k$, then there is a copy of $K_k$ with the $i$th vertex in $V_i$ for each $1\le i\le k$ which is  contained in at least
	\[
		\sum_{j=1}^\ell\left( \prod_{i=1}^k d_{i,k+j}-\delta \right) |V_{k+j}|
	\]
copies of $K_{k+1}$ with the $(k+1)$-th vertex in $\cup_{j=1}^\ell V_{k+j}$.
\end{lemma}

\ignore{We also need the following standard counting lemma, one can see Conlon, Fox and Wigderson \cite[Lemma 2.5]{cfw}, or see Zhao \cite[Theorem~3.27]{Zhao} for a detailed proof.

\begin{lemma}[Conlon, Fox and Wigderson \cite{cfw}]\label{countinglemma}
Suppose that $V_1,\ldots,V_k$ are (not necessarily distinct) subsets of a graph $G$ such that all pairs $(V_i,V_j)$ are $\varepsilon$-regular.
Then the number of labeled copies of $K_k$ whose $i$th vertex is in $V_i$ for all $i$ is at least
\[
\left( \prod_{1 \leq i<j\leq k}d(V_i,V_j) - \varepsilon \binom k2 \right) \prod_{i=1}^k |V_i|.
\]
\end{lemma}
}

We have the following corollary by Conlon, Fox and Wigderson \cite[Corollary 2.6]{cfw}, which counts the monochromatic extensions of cliques.

\begin{corollary}[Conlon, Fox and Wigderson \cite{cfw}]\label{cor:randomclique}
	Let $\varepsilon,\delta \in (0,1)$ and $\varepsilon \leq \delta^3/k^2$. Suppose $U_1,\ldots,U_k$ are (not necessarily distinct) vertex sets in a graph $G$ and all pairs $(U_i,U_j)$ are $\varepsilon$-regular with $\prod_{1 \leq i<j\leq k}d(U_i,U_j) \geq \delta$. Let $Q$ be a randomly chosen copy of $K_k$ with one vertex in each $U_i$ with $1 \leq i \leq k$ and say that a vertex $u$ extends $Q$ if $u$ is adjacent to every vertex of $Q$. Then, for any $u$,
	$\Pr(u \text{ extends }Q) \geq \prod_{i=1}^k d(u,U_i)-4 \delta.$
\end{corollary}

\subsection{General case for $k\ge 3$}
%

Now we give a proof for the upper Ramsey threshold of Theorem~\ref{main} for $k\ge3$.
For any $c>1$ and $k\ge 3$, let $N=c2^kn$ and $p=\frac{1}{c^{1/k}}(1+\gamma),$ where $\gamma>0$ is sufficiently small and $n$ is  sufficiently large. Set
 \[
 p_0= \frac{1}{c^{1/k}}\left(1+\frac{\gamma}{2}\;\right).
  \]
Let $\delta$ and $\varepsilon$ be sufficiently small positive reals such that
\begin{align}\label{param}
\delta=\min\left\{\frac \gamma {4c}, \; \frac{p_0^k}{2^{k+5}}\gamma\right\}, \;\; \text{and} \;\; \varepsilon=\min\left\{\frac{1}{k^2}(\delta p)^k, \; \frac{1}{k^2}(p_0/2)^{\binom k2}\right\}.
\end{align}
Let $\eta>0$ be sufficiently small such that
\begin{align}\label{eta}
\eta=\min\left\{\frac12\delta^k p^k, \;\frac12(p_0/2)^{\binom k2}\right\}.
\end{align}

We begin by applying Lemma \ref{reglem} to the graph $G\in\mathcal{G}(N,p)$ with $\varepsilon$ and $M_0=1/\varepsilon$ to obtain a constant $M=M(\varepsilon)$ such that {\bf a.a.s.} every $2$-coloring of edges of $G\in\mathcal{G}(N,p)$ has an $(\varepsilon, p)$-regular equitable partition $V(G)=\sqcup_{i=1}^m V_i$ where $M_0 \leq m \leq M$ satisfying

\smallskip
(i)	 each part $V_i$ is $(\varepsilon, p)$-regular;

(ii) for each $V_i$, all but at most $\varepsilon m$ parts $V_j$ such that $(V_i,V_j)$ are $(\varepsilon, p)$-regular;

(iii) for any vertex $v\in V$ and for $1\le i\le m$, $\deg_G(v,V_i)\ge p_0|V_i|$;

(iv)  for $1\le i\le m$, $e(V_i)\ge p_0{\binom {|V_i|}2}$;

(v)  for $1\le i<j\le m$, $d_G(V_i,V_j)\ge p_0$.

\medskip

Let $R$ and $B$ be the subgraphs of $G$ induced by all red and blue edges, respectively.
Without loss of generality, we may assume that there are at least $m' \geq m/2$ of the parts, say $V_1,\ldots,V_{m'}$, have internal {\bf red} $p$-density at least $\frac 12$. Let $\Gamma_B$ be the subgraph of the reduced graph $\Gamma$ defined on $\{v_1,\dots,v_m\}$ in which $v_iv_j\in E(\Gamma_B)$ if $(V_i,V_j)$ is $(\varepsilon, p)$-regular and $d_{B,p}(V_i,V_j)\ge \delta$.
Let $\Gamma_B'$ be the subgraph of $\Gamma_B$ induced by the ``red'' vertices $v_i$ for $1 \leq i \leq m'$.

Suppose that, in $\Gamma_B'$, some vertex $v_i$ has at least $(2^{1-k}+2\varepsilon)m'$ non-neighbors. Then, since for $V_i$, there are at most $\varepsilon m \leq 2 \varepsilon m'$ $V_j$'s such that $(V_i,V_j)$ is not $(\varepsilon,p)$-regular, we have that there are at least $2^{1-k}m'$ parts $V_j$ with $1 \leq j \leq m'$ such that $(V_i,V_j)$ is $(\varepsilon,p)$-regular. Let $J$ be the set of all these indices $j$ such that $v_j$ is the non-neighbor of $v_i$ and $(V_i,V_j)$ is $(\varepsilon,p)$-regular. Then $|J|\ge m/2^k$.
Note that
 \[
 d_{B,p}(V_i,V_j)+d_{R,p}(V_i,V_j)=\frac{e_{B}(V_i,V_j)+e_{R}(V_i,V_j)}{p|V_i||V_j|}\ge\frac{p_0}{p},
 \]
 thus if $v_iv_j\not\in E(\Gamma_B)$, then we have $d_{R,p}(V_i,V_j) \geq \frac{p_0}{p}- \delta$ and so the edge density between $V_i$ and $V_j$ satisfies $d_{R}(V_i,V_j) \geq p_0-p\delta$.  Since the red $p$-density is at least $1/2$, from Lemma \ref{con-ct}, there exists a red $K_k$ which is contained in at least
\begin{align*}
\sum_{j\in J}\left((p_0-p\delta)^k-\delta\right)|V_j|&\ge  \left(\left(\frac{1}{c^{1/k}}\left(1+\frac\gamma2\right)-\delta\right)^k-\delta\right)|J|\frac Nm
\\&\ge  \left(\frac{1}{c}\left(1+ck\delta\right)-\delta\right)|J|\frac Nm\ge n
\end{align*}
red $K_{k+1}$ by noting (\ref{param}) that $\delta\le \frac \gamma {4c}$. Thus, we obtain a red $B_n ^{(k)}$ as desired.

Therefore, we may assume that every vertex in $\Gamma_B'$ has degree at least $(1-2^{1-k}- 2\varepsilon)m'$. Since $2^{1-k}+2\varepsilon < \frac{1}{k-1}$ for $k \geq 2$, it follows from Tur\'an's theorem that $\Gamma_B'$ contains a $K_k$ on vertices $v_{i_1},\ldots,v_{i_k}$. Let $W_j=V_{i_j}$ for $1\le j\le k$. Then every pair $(W_i,W_j)$ with $i \leq j$ is $(\varepsilon,p)$-regular and $d_{B,p}(W_i,W_j) \geq \delta$ for $i \neq j$, and each $W_i$ has red $p$-density at least $\frac 12$.

From Lemma \ref{countinglemma} and (\ref{eta}), the number of blue $K_k$'s with the $i$th vertex in $W_i$ is at least
	\begin{align*}
		\left( \prod_{1 \leq i <j \leq k}[p\cdot d_{B,p}(W_i,W_j)] -\eta \right) \prod_{i=1}^k |W_i| &\geq \left( \delta^k p^k- \eta \right) \prod_{i=1}^k |W_i|>0.
	\end{align*}
Similarly, the number of red $K_k$'s in any $W_i$ is at least
	\[
		\left( [p_0\cdot d_{R,p}(W_i)]^{\binom k2} -\eta\right) |W_i|^k
\geq \left( (p_0/2)^{\binom k2}- \eta \right) |W_i|^k>0.
	\]

For any vertex $v$, define
\begin{align*}
d_{B,p}(v,W_i):=\frac{\deg_B(v,W_i)}{p_0|W_i|}.
\end{align*}
 Similarly, we define $d_{R,p}(v,W_i)$. From the assumption that $\deg_G(v,W_i)\ge p_0|W_i|$, we have
\begin{align}\label{r-b-ds}
d_{R,p}(v,W_i)+d_{B,p}(v,W_i)\ge1.
\end{align}	
We may assume that the equality holds for all vertex $v$.

	Now, for any vertex $v$ and for $1\le i \le k$, let
$
x_i(v):=d_{B,p}(v,W_i).
$
Then $d_{R,p}(v,W_i)\ge 1-x_i(v).$
 From  a technical analytic inequality by Conlon \cite[Lemma~8]{Conlon}, we know that
	\[
		\prod_{i=1}^k x_i(v)+ \frac 1k \sum_{i=1}^k (1-x_i(v))^k \geq 2^{1-k}.
	\]

Therefore, we have either $\prod_{i=1}^k x_i(v)\ge 2^{-k}$ or $\frac 1k \sum_{i=1}^k (1-x_i(v))^k\ge 2^{-k}$. There are two cases as follows.

\medskip

{\bf Case 1.} \ $\prod_{i=1}^k x_i(v) \geq 2^{-k}.$

\medskip

	For a given vertex $v$, if we pick $w_i \in W_i$ with $1 \leq i \leq k$ uniformly and independently at random, then the probability that all the edges $(v,w_i)$ are blue is roughly $\prod_{i=1}^k [px_i(v)]$. Together with the regularity of the pairs $(W_i,W_j)$, a random blue $K_k$ spanned by $(W_1,\ldots,W_k)$ will also have probability close to $\prod_{i=1}^k [px_i(v)]$ of being in the blue neighborhood of a random chosen $v$. Indeed, from Corollary \ref{cor:randomclique}, the expected number of blue extensions of a randomly chosen blue $K_k$ spanned by $(W_1,\ldots,W_k)$ is at least
	\begin{align*}
		\sum_{v \in V} \left( \prod_{i=1}^k [p_0\cdot d_{B,p}(v,W_i)]-4 \delta \right)
&=\sum_{v \in V} \left( \prod_{i=1}^k [p_0x_i(v)]-4 \delta \right)
     \geq \left(2^{-k}-\frac{4 \delta}{p_0^k}\right)p_0^kN
  \\&   =\left(2^{-k}-\frac{4 \delta}{p_0^k}\right)\frac1c\left(1+\frac\gamma2\right)^k\cdot c2^kn
		\geq n
	\end{align*}
	by noting $\delta \leq \frac{p_0^k}{2^{k+5}}\gamma$ from (\ref{param}). Therefore, a randomly chosen blue $K_k$ spanned by $(W_1,\ldots,W_k)$ will have at least $n$ blue extensions in expectation, giving us a blue $B_n ^{(k)}$.

\medskip

{\bf Case 2.} \ $\frac 1k \sum_{i=1}^k (1-x_i(v))^k \geq 2^{-k}$.

\medskip

For this case, we have
\[
\frac 1k \sum_{i=1}^k \sum_{v \in V} (1-x_i(v))^k
= \frac 1k \sum_{v \in V}\sum_{i=1}^k  (1-x_i(v))^k \geq 2^{-k}N.
\]
Thus there must exist some $1 \leq i\leq k$ for which $\sum_{v \in V} (1-x_i(v))^k \geq 2^{-k} N$. Similarly, from the regularity of $W_i$, for a random red $K_k$ in $W_i$ and for a random $v \in V$, $v$ will form a red extension of the $K_k$ with probability close to $p^k(1-x_i(v))^{-k}$. Indeed, we can apply Corollary \ref{cor:randomclique} again to obtain that the expected number of extensions of a random red $K_k$ in $W_i$ is at least
	\[
		\sum_{v \in V} \left([p_0(1-x_i(v))]^k-4 \delta\right) \geq (2^{-k}-4 \delta/p_0^k)p_0^kN \geq n,
	\]
yielding a red $B_n ^{(k)}$ as desired. Theorem \ref{main} is proved.\hfill$\Box$

\section{Concluding remarks}\label{clu}

In this paper, we obtain the sharp Ramsey threshold for the book graph $B_n^{(k)}$. In particular, for every fixed integer $k\ge 2$ and for any real $c>1$, let $N=c2^k n$. Then for any real $\gamma>0$,
\[
\lim_{n\to \infty} \Pr(G(N,p)\to B_n^{(k)})= \left\{  \begin{array}{cl}
0  & \mbox{if $p\le\frac{1}{c^{1/k}}(1-\gamma)$,} \\  1 & \mbox{if $p\ge\frac{1}{c^{1/k}}(1+\gamma)$}. \end{array}      \right.
\]
Note that $N=c2^k n=(c+o(1))r(B_n^{(k)},B_n^{(k)})$. It would be interesting to determine the sharp Ramsey threshold for the book graph $B_n^{(k)}$ if $N=\omega(n)\cdot r(B_n^{(k)},B_n^{(k)})$, where $\omega(n)$ tends to infinity as $n\to\infty$. As a special case, it would be interesting to determine the sharp Ramsey threshold for the book graph when $N=\Theta(n^2)$.


\end{spacing}


\begin{thebibliography}{99}

\bibitem{alo-spe}
N. Alon and J. Spencer, The Probabilistic Method,
Wiley-Interscience, New York, 1992.



\bibitem{beck-83}
J. Beck, On size Ramsey number of paths, trees, and circuits, I, {\em J. Graph Theory} 7 (1983), 115--129.

\bibitem{beck-2}
J. Beck, On Size Ramsey Number of Paths, Trees and Circuits II, in: Mathematics of Ramsey Theory (Nes\u{e}t\u{r}l and R\"{o}dl eds.), 34--45, Springer-Verlag, 1990.

\bibitem{Bollobas}
B. Bollob\'as, Random Graphs, Cambridge University Press, 2001.



\bibitem{cgms}
M. Campos, S. Griffiths, R. Morris  and J. Sahasrabudhe, An exponential improvement for diagonal Ramsey, arXiv: 2303.09521v1, 2023.

\bibitem{cl}
X. Chen and Q. Lin, New upper bounds for Ramsey numbers of books, {\em European J. Combin.} 115 (2024), Paper No. 103785, 9 pp.


\bibitem{cly}
X. Chen, Q. Lin and C. You, Ramsey numbers of large books, {\em J. Graph Theory} 101 (2022), no. 1, 124--133.

\bibitem{cher}
H. Chernoff, A measure of the asymptotic efficiency for tests of a
hypothesis based on the sum of observations, {\em Ann. Math.
Statistics} 23 (1952), 493--507.


\bibitem{c-r-s-t}
C. Chvat\'al, V. R\"odl, E. Szemer\'edi and W. Trotter, The Ramsey number of a graph with bounded maximum degree, {\em J. Combin. Theory Ser. B} 34 (1983), 239--243.


\bibitem{Conlon}
D.~Conlon, The {Ramsey} number of books, \emph{{A}dv. {C}ombin.}  (2019), Paper
  No. 3, 12 pp.

\bibitem{cfw}
D. Conlon, J. Fox and Y. Wigderson, Ramsey numbers of books and quasirandomness, {\em Combinatorica} 42 (2022), no. 3, 309--363.

\bibitem{cfw-2}
D. Conlon, J. Fox and Y. Wigderson, Off-diagonal book Ramsey numbers, {\em Combin. Probab. Comput.} 32 (2023), 516--545.

\bibitem{cfw-s}
D. Conlon, J. Fox and Y. Wigderson,
Three early problems on size Ramsey numbers, to appear in Combinatorica.

\bibitem{d-p}
A. Dudek and P. Pra{\l}at, An Alternative proof of the linearity of the size-Ramsey number of paths, {\em Combin. Probab. Comput.} 24 (2015), 551--555.


\bibitem{efrs-size}
P. Erd\H{o}s, R. Faudree, C. Rousseau and R. Schelp, The size Ramsey numbers, {\em Period. Math. Hungar.} 9 (1978), 145--161.


\bibitem{erd-ren}
P. Erd\H{o}s and A. R\'enyi, On the evolution of random graphs, {\em Publ. Math. Inst. Hungar. Acad. Sci.} 5 (1960), 17--61.


\bibitem{fly}
C. Fan, Q. Lin and Y. Yan, On a conjecture of Conlon, Fox, and Wigderson, Combin. Probab. Comput. (2024), doi:10.1017/S0963548324000026.

\bibitem{f-r}
P. Frankl and V. R\"{o}dl, Large triangle-free subgraphs in graphs without $K_4$, {\em Graphs Combin.} 2 (1986), 135--144.

\bibitem{f-k}
E. Friedgut and M. Krivelevich, Sharp thresholds for certain Ramsey properties of random graphs,
{\em Random Structures Algorithms} 17 (2000), no. 1, 1--19.

\bibitem{f99}
E. Friedgut, Sharp thresholds of graph properties, and the $k$-sat problem, {\em J. Amer. Math. Soc.} 12 (1999) no. 4, 1017--1054.

\bibitem{frrt}
E. Friedgut, V. R\"{o}dl, A. Ruci\'{n}ski and P. Tetali, A sharp threshold for random graphs with a monochromatic triangle in every edge coloring, {\em Mem. Amer. Math. Soc.} 179 (2006), no. 845, vi+66 pp.

\bibitem{g-g}
L. Gerencs\'er and A. Gy\'arf\'as, On Ramsey-type problems, {\em Ann. Univ. E\"{o}tv\"{o}s Sect. Math.} 10 (1967), 167--170.

\bibitem{jlr}
S. Janson, T. {\L}uczak and A. Ruci\'nski, Random Graphs, Wiley-Interscience, New York, 2000.

\bibitem{k-r1}
Y. Kohayakawa, Szemer\'{e}di's regularity lemma for sparse graphs, Foundations of computational mathematics (Rio de Janeiro, 1997), Springer, Berlin, 1997, 216--230.





\bibitem{k-r-2}
Y. Kohayakawa and V. R\"{o}dl, Szemer\'{e}di's regularity lemma and quasi-randomness, Recent advances in algorithms and combinatorics, CMS Books Math./Ouvrages Math. SMC, vol. 11, Springer, New York,
2003, 289--351.

\bibitem{k-r-s-s}
Y. Kohayakawa, V. R\"{o}dl, M. Schacht and E. Szemer\'edi,
Sparse partition universal graphs for graphs of bounded degree, {\em Adv. Math.} 226 (2011), 5041--5065.

 \bibitem{kss}
J. Koml\'{o}s, A. Shokoufandeh, M. Simonovits and E. Szemer\'{e}di, The regularity lemma and its
applications in graph theory, Theoretical aspects of computer science (Tehran, 2000), Lecture
Notes in Comput. Sci., vol. 2292, Springer, Berlin, 2002, pp. 84--112.

\bibitem{ko-sim}
J. Koml\'{o}s and M. Simonovits, Szemer\'{e}di's regularity lemma and its applications to graph theory.
{\em Combinatorics, Paul Erd\H{o}s is eighty, Vol. 2 (Keszthely, 1993)}, 295--352, Bolyai Soc. Math. Stud., 2, {\em J\'{a}nos Bolyai Math. Soc., Budapest}, 1996.

\bibitem{let}
S. Letzter, Path Ramsey number for random graphs, {\em Combin. Probab. Comput.} 25 (2016), 612--622.

\bibitem{li-l}
Y. Li and Q. Lin, Elementary methods of graph Ramsey theory, Springer, 2022.


\bibitem{ltz}
Y. Li, X. Tang and W. Zang, Ramsey functions involving $K_{m,n}$ with $n$ large, {\em Discrete Math.} 300 (2005) no. 1-3, 120--128.

\bibitem{l-r-v}
T. {\L}uczak, A. Ruci\'nski and B. Voigt, Ramsey properties of random graphs, {\em J. Combin. Theory Ser. B}  56 (1992), 55--68.



\bibitem{n-s}
R. Nenadov and A. Steger, A short proof of the random Ramsey theorem, {\em Combin. Probab. Comput.} 25 (2016), 130--144.

\bibitem{n1} V. Nikiforov and C. Rousseau, A note on Ramsey numbers for books, {\em J. Graph Theory} 49 (2005), 168--176.

\bibitem{n2} V. Nikiforov and C. Rousseau, Book Ramsey numbers I, {\em Random Structures Algorithms} 27 (2005),  379--400.

\bibitem{n3} V. Nikiforov, C. Rousseau and R. Schelp, Book Ramsey numbers and quasi-randomness, {\em Combin. Probab. Comput.} 14 (2005), 851--860.


 \bibitem{ram}
F.~P.~Ramsey, On a problem of formal logic,  {\em Proc.~Lond. Math.~Soc.} 30 (1929), 264--286.



\bibitem{r-r93}
V. R\"{o}dl and A.  Ruci\'nski, Lower bounds on probability thresholds
for Ramsey properties, {\em Combinatorics, Paul Erd\H{o}s is Eighty} (Vol.1),
Keszthely (Hungary), Bolyai Soc. Math. Studies, 1993, pp.317--346.

\bibitem{r-r-1}
V. R\"{o}dl and A.  Ruci\'nski, Random graphs with monochromatic triangles in every edge coloring, {\em Random Structures Algorithms} 5 (1994), 253--270.

\bibitem{r-r}
V. R\"{o}dl and A. Ruci\'nski, Threshold functions for Ramsey properties, {\em J. Amer. Math. Soc.} 8 (1995), 917--942.

 \bibitem{rs}
V. R\"{o}dl  and M. Schacht, Regularity lemmas for graphs, in: {\em Fete of Combinatorics
and Computer Science}, Bolyai Soc. Math. Stud. 20, 2010, 287--325.


\bibitem{r-s}
V. R\"{o}dl and E. Szemer\'edi, On size Ramsey numbers of graphs with bounded degree, {\em Combinatorica} 20 (2000), 257--262.

\bibitem{rou-sh}
C. Rousseau and J. Sheehan, On Ramsey numbers for books, {\em J. Graph Theory} 2 (1978) 77-87.

\bibitem{szemeredi}
E. Szemer\'edi, On sets of integers containing no $k$ elements in arithmetic progression, {\em Acta Arith.} 27 (1975), 199--245.

\bibitem{szemeredi1}
E. Szemer\'edi, Regular partitions of graphs, in Probl\'emes combinatories et th\'eorie des graphes, Colloq. Internat., CNRS, 260, Paris, 1978, 399--401.


 \bibitem{tho}
A.~Thomason, On finite Ramsey numbers, {\em European J.~Combin.} {3} (1982), 263--273.



\bibitem{Zhao}
Y.~Zhao, Graph theory and additive combinatorics: Notes for {MIT} 18.217, 2019.
 http://yufeizhao.com/gtac/gtac.pdf.



\end{thebibliography}
\end{document}